\documentclass[12pt]{article}
\usepackage{latexsym}

\newcommand{\duk}{\noindent {\bf Proof. }}
\newcommand{\kduk}{\hfill $\Box$\bigskip}

\newcommand{\N}{\mathbf{N}}

\newcommand{\ep}{\varepsilon}
\newcommand{\al}{\alpha}

\newcommand{\st}{\mathrm{al}}
\newcommand{\ind}{\mathrm{Ind}}
\newcommand{\po}{\mathrm{Seq}}

\newtheorem{veta}{Theorem}[section]

\newtheorem{lema}[veta]{Lemma}

\newtheorem{prop}[veta]{Proposition}

\def\cla#1#2#3#4#5#6{%autor, nazev, cas., rocnik, rok, stranky
  {\sc #1, }#2, {\it #3, }{\bf #4 }(#5), #6.}

\def\kni#1#2#3#4#5{%autor, nazev, nakladatel, sidlo, rok
  {\sc #1, }{\it #2, }#3, #4, #5.}

\def\vsbo#1#2#3#4#5#6#7#8{%autor, nazev, stranky, editor, nazev, nakladatel, 
  {\sc #1, }#2. In: {#4 (ed.), } {\it #5, } #6, #7,  #8; pp. #3.}

\begin{document}

\author{Martin Klazar\thanks{Department of Applied Mathematics (KAM) and Institute for Theoretical 
Computer Science (ITI), Charles University, Malostransk\'e n\'am\v est\'\i\ 25, 118 00 Praha, 
Czech Republic. ITI is supported by the project LN00A056 of the 
Ministry of Education of the Czech Republic. E-mail: {\tt klazar@kam.mff.cuni.cz}}}
\title{On the least exponential growth admitting uncountably many closed permutation classes}
\date{}

\maketitle
\begin{abstract}
We show that the least exponential growth of counting functions which admits uncountably many closed 
permutation classes lies between $2^n$ and $(2.33529\dots)^n$.
\end{abstract}

\section{Introduction}

Let $S_n$ be the set of $n!$ permutations of $[n]=\{1,2,\dots,n\}$, $S=\bigcup_{n=0}^{\infty}S_n$ 
be the set of all finite permutations, and $\prec$ be the usual containment of permutations (defined below).
It is well-known that the partial ordering $(S,\prec)$ has infinite antichains, see  
\cite{lave}, \cite{prat}, \cite{spie_bona}, and \cite{tarj}. Equivalently, $(S,\prec)$ has uncountably many 
lower order ideals $X\subset S$; these are called {\em closed permutation classes} or, for short, CPC's. 
In this article we 
want to localize the least exponential growth of the counting function $n\mapsto |X\cap S_n|$ 
which admits uncountably many CPC's $X$. 

More precisely, if
$$
K_{\al}=\{X:\ \mbox{$X$ is a CPC such that $|X\cap S_n|<\al^n$ for all $n>n_0$}\},
$$
what can be said about the number
$$
\kappa=\inf\{\al>1:\ \mbox{the set $K_{\al}$ is uncountable}\}.
$$
We prove the following bounds. 
\begin{veta}\label{hlavni}
Let $\kappa$ determine the least exponential growth of uncountably many CPC's, as defined above. Then
$$
2\le\kappa\le 2.33529\dots
$$
where the upper bound is the only real root of $x^5-x^4-2x^3-2x^2-x-1$. 
\end{veta}
When the base $\al$ in $\al^n$  is increased, the ``phase transition" from countably to 
uncountably many CPC's with growth $<\al^n$, $n>n_0$, occurs somewhere in the interval $[2,2.33529\dots]$. 
It would be interesting to narrow it or to determine $\kappa$ exactly. 

In the proof of Theorem~\ref{hlavni} we build on previously obtained results. In Kaiser and 
Klazar \cite[Theorem 3.8]{kais_klaz} we have proved that the exponential growths of CPC's $X$
such that $|X\cap S_n|<2^{n-1}$ for at least one $n$ form a discrete hierarchy $\al_i^n$, $i=2,3,4,\dots$, 
where $\al_2=1.61803\dots<\al_3<\al_4<\dots<2$, 
$\al_i\uparrow 2$, and $\al_i$ is the largest positive real root of $x^i-x^{i-1}-\cdots-1$. It follows from
the proof, with some additional arguments from the wqo theory, that the structure of the corresponding 
CPC's is so restricted that each set $K_{2-\ep}$ must be countable. In Section 2 we give a proof of this 
fact. On the other hand, Spielman and B\'ona \cite{spie_bona} constructed an infinite antichain $(R,\prec)$ 
such that $123\not\prec\pi$ for every $\pi\in R$. Thus, denoting $S(123)$ the set of 
$123$-avoiding permutations, 
there are uncountably many CPC's $X$ with $X\subset S(123)$. Since 
$|S(123)\cap S_n|=\frac{1}{n+1}{2n\choose n}$ (Rogers \cite{roge}, Simion and Schmidt \cite{simi_schm}, 
$\dots$), we obtain the bound $\kappa\le 4$. The enumeration of $S(123,3214)$, due to West \cite{west96}, 
and the infinite antichain $U$ due to Atkinson, Murphy and Ru\v skuc \cite{atki02} give the improvement 
$\kappa\le 2.61803\dots$. In Section 3 we lower this further to the upper bound 
in Theorem~\ref{hlavni}.

Closed permutation classes and permutation avoidance (containment) are related to 
computer science mainly via sorting problems. The set of permutation $\pi$ which, when inputed to 
some sorting device,
can be sorted to the identical permutation, is often a CPC. Indeed, this was the very first 
motivation to introduce $\prec$ in the works of Pratt \cite{prat} and Tarjan \cite{tarj}. Recent works 
on closed permutation classes and permutation containment with motivation in computer science 
(sorting, complexity of recognizing $\prec$) 
are, for example, Ahal and Rabinovich \cite{ahal}, Albert et al. \cite{albe01}, Atkinson \cite{atki98}, 
Atkinson, Murphy and Ru\v skuc \cite{atki02a, atki03}, and Bose, Buss and Lubiw \cite{bose}.

Now we review the definition of $\prec$ and basic facts on CPC's.
Further definitions will be given throughout next two sections.  

For $\pi\in S_n$, $n$ is the {\em length} of $\pi$ and we define $|\pi|=n$. 
For $A,B\subset\N=\{1,2,\dots\}$ the notation $A<B$ means that $a<b$ for every $a\in A$ and $b\in B$. Interval 
$\{a,a+1,a+2,\dots,b\}$, where $a,b\in\N$, is denoted $[a,b]$. Instead of $[1,n]$ we write $[n]$.
Two $m$-term sequences $a_1a_2\dots a_m$ and $b_1b_2\dots b_m$ over $\N$ are {\em order-isomorphic} if 
$b_k<b_l\Leftrightarrow a_k<a_l$ for all $k,l\in[m]$. A permutation $\pi$ is {\em contained} in another 
permutation $\rho$, written $\pi\prec\rho$, if $\rho$ (as a sequence) has a subsequence that is 
order-isomorphic to $\pi$; in the opposite case $\rho$ is $\pi$-{\em avoiding}. Visually, the graph of $\pi$ 
(as a discrete function) can be obtained
from that of $\rho$ by omitting points. If $\pi\in S_n$ and $A\subset[n]$, the 
{\em restriction} $\pi|A$ is the 
permutation order-isomorphic to the corresponding subsequence of $\pi$.  
For $X\subset S$, $M(X)$ is the set of all $\prec$-minimal 
permutations not in $X$, and $S(X)$ is the set of all permutations not containing any member of $X$. We define 
$S_n(X)=S(X)\cap S_n$. 
For finite $X=\{\pi_1,\dots,\pi_r\}$ we write $S(\pi_1,\dots,\pi_r)$ and $S_n(\pi_1,\dots,\pi_r)$ instead of 
$S(\{\pi_1,\dots,\pi_r\})$ and $S_n(\{\pi_1,\dots,\pi_r\})$.
Clearly, each proper restriction of each $\pi\in M(X)$ lies in $X$.
A set $X\subset S$ is a CPC (closed permutation class) if $\pi\prec\sigma\in X$ implies $\pi\in X$. 
Each $S(X)$ is a CPC and for each CPC $X$ we have $X=S(M(X))$. Each $M(X)$ is an antichain 
(its elements are mutually incomparable by $\prec$) and for each antichain $X\subset S$ we have $X=M(S(X))$.  
Thus the mapping $X\mapsto M(X)$, with the inverse 
$X\mapsto S(X)$, is a bijection between the set of all CPC's and the set of all antichains of permutations. 

\section{The lower bound of Theorem~\ref{hlavni}}

In this section we mostly follow the notation of \cite{kais_klaz}. 
A permutation $\sigma$ is {\em alternating} if $\sigma(\{1,3,5,\dots\})>\sigma(\{2,4,6,\dots\})$. 
For $\pi\in S$ we let $\st(\pi)$ be the maximum 
length of an alternating permutation $\sigma$ such that $\sigma\prec\pi$ or $\sigma\prec\pi^{-1}$. For a set 
of permutations $X$ we denote $\st(X)=\max\{\st(\pi):\ \pi\in X\}$. 

\begin{lema}\label{stridavost}
If $X$ is a CPC with $\st(X)=\infty$, then $|X\cap S_n|\ge 2^{n-1}$ for every $n\in\N$. 
\end{lema}
\duk
We suppose that $X$ contains arbitrarily long alternating permutations; the other case with inverses is 
treated similarly.
Using the closeness of $X$ and the pigeonhole principle, we deduce that either for every $n\in\N$ there is an 
alternating $\pi\in X\cap S_n$ such that $\pi(1)<\pi(i)$ for every odd $i\in[2,n]$ or 
for every odd $n\in\N$ there is an alternating $\pi\in X\cap S_n$ such that $\pi(n)<\pi(i)$ for every odd 
$i\in[n-1]$. We assume 
that the former case occurs, the latter one is similar. It follows that for every $n\in\N$ and every subset 
$A\subset[2,n]$ there is a permutation $\pi_A\in X\cap S_n$ such that $\pi_A(i)<\pi_A(1)\Leftrightarrow 
i\in A$. For 
distinct subsets $A$ we get distinct permutations $\pi_A$ and $|X\cap S_n|\ge 2^{n-1}$.
\kduk

If $\sigma\in S_n$ and $\tau\in S_m$, then $\pi=\sigma\oplus\tau\in S_{n+m}$ is the permutation defined 
by $\pi(i)=\sigma(i)$ for $i\in[n]$ and $\pi(i)=n+\tau(i-n)$ for $i\in[n+1,n+m]$. Similarly, 
$\pi=\sigma\ominus\tau$ is defined by $\pi(i)=m+\sigma(i)$ for $i\in[n]$ and $\pi(i)=\tau(i-n)$ for 
$i\in[n+1,n+m]$. Note that if  $\sigma'\prec\sigma$ and $\tau'\prec\tau$, then 
$\sigma'\oplus\tau'\prec\sigma\oplus\tau$; similarly for $\ominus$.
If $\pi\in S$ has no decomposition $\pi=\sigma\oplus\tau$ for any nonempty $\sigma$ and 
$\tau$, we say that $\pi$ is {\em up-indecomposable}. The subset of up-indecomposable permutations in 
$S_k$ is denoted $\ind_k^+$. Each $\pi\in S$ has a unique {\em up-decomposition}
$\pi=\sigma_1\oplus\sigma_2\oplus\dots\oplus\sigma_k$ where each $\sigma_i$ is up-indecomposable; $\sigma_i$'s 
are called {\em up-blocks}. 
The maximum size of an up-block in the up-decomposition of $\pi$ is denoted $h^+(\pi)$. 
For the operation $\ominus$, the down-(in)decomposability, sets $\ind_k^-$, down-decompositions, down-blocks, 
and function $h^-(\cdot)$ are defined in an analogous way. 

The proof of the next lemma is left to the reader as an exercise (or see \cite[Lemma 3.7]{kais_klaz}).

\begin{lema}\label{zmenseni}
For every $\pi\in \ind_n^+$, $n>1$, there is a $\sigma\in \ind_{n-1}^+$ such that $\sigma\prec\pi$. The 
same holds for down-indecomposable permutations.
\end{lema}

\begin{lema}\label{vlacek}
If $X$ is a CPC with the property that for every $k\in\N$ there is a permutation $\sigma\in\ind_k^+$ 
such that $\sigma\oplus\sigma\oplus\dots\oplus\sigma\in X$ ($k$ summands), then $|X\cap S_n|\ge 2^{n-1}$ for 
every $n\in\N$.  An analogous result holds for down-decompositions.
\end{lema}
\duk
Using the assumption and Lemma~\ref{zmenseni}, we obtain that for every $n\in\N$ there is a set 
$\Sigma=\{\sigma_1,\sigma_2,\dots,\sigma_n\}$ such that $\sigma_i\in\ind_i^+$ and every permutation 
of the form $\pi=\rho_1\oplus\rho_2\oplus\dots\oplus\rho_r$, where $\rho_i\in\Sigma$ and $r\le n$, 
is in $X$. Since the up-decomposition uniquely determines $\pi$, there are exactly $2^{n-1}$ such 
permutations $\pi$ in $X\cap S_n$ (as compositions of $n$) and $|X\cap S_n|\ge 2^{n-1}$.
\kduk

Let $H^+_k=\{\pi\in S:\ h^+(\pi)<k\}$ and similarly for $H^-_k$.
For $k\in\N$ and $\pi\in S_n$, we let $s_k(\pi)$ be the number $r$ of intervals $I_1<I_2<\dots<I_r$ 
in this unique decomposition of $[n]$: $I_1$ is the longest initial interval in $[n]$ such that 
$\pi|I_1\in H^+_k\cup H^-_k$, $I_2$ is the longest following interval such that 
$\pi|I_2\in H^+_k\cup H^-_k$ and so on. We call $I_1<I_2<\dots<I_r$ the $k$-{\em decomposition} of $\pi$. 
Note that each restriction $\pi|I_i$ has up-decomposition or down-decomposition composed of blocks of lengths at most 
$k-1$ and that each restriction $\pi|I_i\cup I_{i+1}$ contains both an element from $\ind_k^+$ and 
an element from $\ind_k^-$. For $k\in\N$ and $X$ a set of permutations we define 
$s_k(X)=\max\{s_k(\pi):\ \pi\in X\}$. We let $s_1(\pi)=s_1(X)=\infty$ for every permutation $\pi$ and set $X$. 

\begin{prop}\label{struktura}
If $X$ is a CPC such that $|X\cap S_n|<2^{n-1}$ for some $n\in\N$, then $\st(X)<\infty$ and, for some 
$k\in\N$, $s_k(X)<\infty$.  
\end{prop}
\duk
If $\st(X)=\infty$, we have $|X\cap S_n|\ge 2^{n-1}$ for all $n\in\N$ by Lemma~\ref{stridavost}, which is 
a contradiction. Suppose that $s_k(X)=\infty$ for every $k\in\N$. By the remark after the definition of 
$s_k(\cdot)$, the pigeonhole principle and the closeness of $X$, for every $k\ge 2$ there are permutations 
$\sigma_k\in\ind_k^+$, $\tau_k\in\ind_k^-$ and 
$\pi_{k}\in X\cap S_r$, $k^2\le r\le 2k^2$, with the property that $[r]$ can be decomposed into $k$ intervals 
$I_{k,1}<I_{k,2}<\dots<I_{k,k}$, 
$k\le |I_{k,i}|\le 2k$, so that each of the $k$ restrictions $\pi_{k}|I_{k,i}$ contains both $\sigma_k$ and 
$\tau_k$. For $k\in\N$ and $1\le i\le k$, we consider the interval
$$
J_{k,i}=[\min\pi_{k}(I_{k,i}),\max\pi_{k}(I_{k,i})].
$$
Using the Ramsey theorem and Lemma~\ref{zmenseni}, we may assume that either for every $k\in\N$ the $k$ 
intervals $J_{k,1},\dots,J_{k,k}$ intersect each other or for every $k\in\N$ these $k$ intervals are 
mutually disjoint. In the former case, they must always have one point in common, and it follows that 
$\st(X)=\infty$. 
We have again the contradiction by Lemma~\ref{stridavost}. In the latter case, using again Ramsey theorem 
(or Erd\H os--Szekeres theorem) and Lemma~\ref{zmenseni}, we may assume that either for every 
$k\in\N$ we have 
$J_{k,1}<J_{k,2}<\dots<J_{k,k}$ or for every $k\in\N$ we have $J_{k,1}>J_{k,2}>\dots>J_{k,k}$. 
Then for every 
$k\in\N$ we have $\sigma_k\oplus\sigma_k\oplus\dots\oplus\sigma_k\in X$ ($k$ summands) or for every 
$k\in\N$ we have $\tau_k\ominus\tau_k\ominus\dots\ominus\tau_k\in X$ ($k$ summands). By 
Lemma~\ref{vlacek}, we get the contradiction that $|X\cap S_n|\ge 2^{n-1}$ for all $n\in\N$.
\kduk

Every bijection $f:X\to Y$, where $X=\{x_1<x_2<\dots<x_n\}$ and $Y=\{y_1<y_2<\dots<y_n\}$ are subsets 
of $\N$, defines a unique $\pi\in S_n$ {\em order-isomorphic} to $f$: $\pi(i)=j\Leftrightarrow f(x_i)=y_j$. 
An interval in $X$ is a subset of the form $\{x_i,x_{i+1},\dots,x_{j}\}$, $1\le i\le j\le n$.

\begin{lema}\label{zjemneni}
Let $X,Y\subset\N$ be two $n$-element subsets, $f:X\to Y$ be a bijection, and $\pi\in S_n$ be 
order-isomorphic to $f$. Suppose $\pi\in H^+_k\cup H^-_k$. Then every interval partition $J_1<J_2<\dots<J_r$ 
of $X$ can be refined by an interval partition $I_1<I_2<\dots<I_s$ such that $s\le r+(k-1)(r-1)$ and 
each image $f(I_i)$ is an interval in $Y$. Similarly, every partition of $Y$ in $r$ intervals can be refined 
by a partition in at most $r+(k-1)(r-1)$ intervals which under $f^{-1}$ map to intervals in $X$.
\end{lema}
\duk
It suffices to prove only the first part because $\pi\in H^+_k\cup H^-_k$ implies that 
$\pi^{-1}\in H^+_k\cup H^-_k$. Without loss of generality we can assume that $X=Y=[n]$ and $f=\pi$.
Let $\pi\in S_n\cap H^+_k$ (the case with $H^-_k$ is similar) and $J_1<J_2<\dots<J_r$ be an interval 
partition of $[n]$. We call an up-block in the up-decomposition
$\pi=\sigma_1\oplus\sigma_2\oplus\dots\oplus\sigma_t$ intact if its domain lies completely in some $J_i$ and 
we call it split otherwise. Clearly, there are at most $r$ maximal runs of intact up-blocks and 
at most $r-1$ split up-blocks. We partition $[n]$ in the intervals $I_1<I_2<\dots<I_s$ so that each $I_i$ is either 
the domain of a maximal run or a singleton in the domain of a split up-block. Since $|\sigma_i|<k$ for each $i$, 
we have $s\le r+(k-1)(r-1)$. This is a refinement of the original interval partition and $\pi(I_i)$ is 
an interval for every $i$. 
\kduk

We will need a continuity property of the functions $\st(\cdot)$ and $s_k(\cdot)$.

\begin{lema}\label{spojitost}
Let $\sigma\in S_n$, $\tau\in S_{n+1}$, and $\sigma\prec\tau$. Then $\st(\tau)\le \st(\sigma)+2$ and, 
for every $k\in\N$, $s_k(\tau)\le s_k(\sigma)+2$. 
\end{lema}
\duk
Let $\rho\in S$ be alternating, $|\rho|=\st(\tau)$, and $\rho\prec\tau$ (the case $\rho\prec\tau^{-1}$
is similar).
The permutation $\sigma$ arises by deleting one point from the graph of $\tau$. If this point does not 
lie in the embedding of $\rho$ in $\tau$, we have $\rho\prec\sigma$ and $\st(\sigma)\ge \st(\tau)$. If it does, 
we can delete one more point from the graph of $\rho$ so that the resulting $\rho'$ 
is alternating. But $\rho'\prec\sigma$ and $|\rho'|=|\rho|-2$, so $\st(\sigma)\ge \st(\tau)-2$.

Let $k\ge 2$ be given, $\pi\in S_n$ be arbitrary, and $I_1<I_2<\dots<I_s$ be any decomposition of 
$[n]$ into $s$ intervals satisfying, for every $i=1,\dots,s$, 
$\pi|I_i\in H^+_k\cup H^-_k$; this can be called a weak $k$-decomposition of $\pi$. 
We claim that $s_k(\pi)\le s$. This follows from the 
observation that each interval of the $k$-decomposition of $\pi$ must contain 
the last element of some $I_i$. Now $\tau$ arises by inserting a new point $p$ in the graph of $\sigma$.
The domain $\{p_0\}$ of $p$ is inserted in an interval $J_j$ of the $k$-decomposition $J_1<J_2<\dots<J_r$ 
of $\sigma$  and splits it 
into three intervals $J_j'$, $\{p_0\}$, and $J_j''$ ($J_j'$ or $J_j''$ may be empty). Replacing $J_j$ by 
$J_j'$, $\{p_0\}$, and $J_j''$, we get a weak $k$-decomposition of $\tau$ with at most $r+2$ intervals. 
Thus $s_k(\tau)\le r+2=s_k(\sigma)+2$.
\kduk

Recall that a partial ordering $(Q,\le_Q)$ is a {\em well partial ordering}, briefly wpo, if it has no 
infinite strictly descending chains and no infinite antichains. The first condition is in $(S,\prec)$ satisfied 
but the second one is not and therefore $(S,\prec)$ is not a wpo. Let $(Q,\le_Q)$ be a partial 
ordering. The set $\po(Q)$ of all finite tuples 
$(q_1,q_2,\dots,q_m)$ of elements from $Q$ is partially ordered by the derived Higman ordering $\le_H$: 
$(q_1,q_2,\dots,q_m)\le_H(r_1,r_2,\dots,r_n)\Leftrightarrow$ there is an increasing mapping $f:[m]\to[n]$ such that 
$q_i\le_Q r_{f(i)}$ for every $i\in[m]$. For the proof of the following theorem see Higman \cite{higm} 
or Nash-Williams \cite{nash}.

\begin{veta}[Higman, 1952]\label{higman}
If $(Q,\le_Q)$ is a wpo then $(\po(Q),\le_H)$ is a wpo as well. 
\end{veta}

If $\sigma\in S_m$ and $\tau_i\in S_{n_i}$, $i=1,\dots,m$, the permutation 
$\pi=\sigma[\tau_1,\dots,\tau_m]\in S_{n_1+\cdots+n_m}$ is defined, for $i\in[n_1+\cdots+n_m]$ and setting 
$k=\max(\{j:\ n_1+\cdots+n_j<i\}\cup\{0\})$ and $n_0=0$, by 
$$
\pi(i)=n_0+n_1+\cdots+n_k+\tau_{k+1}(i-n_0-n_1-\cdots-n_k).
$$
Visually, for $i=1,\dots,m$ the $i$-th point (counted from the left) in the graph of $\sigma$ is replaced by a 
downsized copy of the graph of $\tau_i$; the copies are small enough not to interfere 
horizontally and vertically each with the other.  
This operation generalizes $\oplus$ and $\ominus$: $\sigma\oplus\tau=12[\sigma,\tau]$ and 
$\sigma\ominus\tau=21[\sigma,\tau]$. If $\tau_i'\prec\tau_i$, $i=1,\dots,m$, then 
$\sigma[\tau_1',\dots,\tau_m']\prec\sigma[\tau_1,\dots,\tau_m]$.
If $P$ and $Q$ are sets of permutations, we define 
$$
P[Q]=\{\pi[\sigma_1,\dots,\sigma_m]:\ m\in\N,\pi\in P\cap S_m, \sigma_i\in Q\}.
$$ 

The next lemma is an immediate consequence of Higman's theorem or of the easier result that the Cartesian 
product of two wpo's also is a wpo. 

\begin{lema}\label{easy}
Let $P$ and $Q$ be sets of permutations such that $P$ is finite and $(Q,\prec)$ is a wpo. Then
$(P[Q],\prec)$ is a wpo.
\end{lema}

\noindent
Let $\pi\in S_n$ and $J_1<J_2<\dots<J_r$ be an interval partition of $[n]$. Observe that if each image 
$\pi(J_i)$ is also an interval, then there is a permutations $\sigma\in S_r$ such that 
$\pi=\sigma[\pi|J_1,\dots,\pi|J_r]$. 

\begin{lema}\label{substituce}
For every fixed $k,K\in\N$ there is a finite set of permutations $P$ such that 
$$
\{\pi\in S:\ \st(\pi)<K\;\&\;s_k(\pi)<K\}\subset P[H^+_k\cup H^-_k].
$$
\end{lema}
\duk
We show that 
$$
P=S_1\cup S_2\cup\dots\cup S_{kK^*}
$$ 
works where $K^*=(K-1){K\choose 2}+1$. Let $\pi\in S_n$ satisfy $\st(\pi)<K$ and $s_k(\pi)<K$. 
Since $s_k(\pi)<K$, $[n]$ can be partitioned 
in $r$ intervals $J_1<J_2<\dots<J_r$, $r<K$, so that always $\pi|J_i\in H^+_k\cup H^-_k$ 
(we will not need the other property of $k$-decomposition of $\pi$). We show that $[n]$ can be partitioned 
in at most $kK^*$ intervals so that their images under $\pi^{-1}$ are intervals refining
$J_1<J_2<\dots<J_r$. Then we are done because $\pi|I\in H^+_k\cup H^-_k$ for every interval (in fact, every
subset) $I\subset J_i$. 

We consider two words $u$ and $u'$ over $[K]$. The word $u=a_1a_2\dots a_n$ is defined by 
$a_i=j\Leftrightarrow\pi^{-1}(i)\in J_j$ and $u'$ arises from $u$ by contracting each maximal run of one letter 
in one element. For example, if $u=2221331111$ then $u'=2131$. Let $l$ be the length of $u'$ which is also the 
number of maximal runs in $u$.
Clearly, $u'$ has no two consecutive identical letters. Since $\st(\pi)<K$, $u$ 
and $u'$ have no alternating subsequence $\dots a\dots b\dots a\dots b\dots$, $a\ne b$, of length $K+1$. 
A pigeonhole argument implies that $l\le K^*=(K-1){K\choose 2}+1$. 

We partition $[n]$ in $l$ intervals $L_1<L_2<\dots<L_l$ according to the maximal runs in $u$. Each 
$\pi^{-1}(L_i)$ is a subset of some $J_j$ but in general is not an interval. Let $j\in[r]$ and $M_j\subset[n]$ 
be the union of $i_j$ intervals $L_i$ corresponding to all $i_j$ maximal runs of $j$ in $u$; 
$\pi^{-1}(M_j)=J_j$. Applying Lemma~\ref{zjemneni} to the restricted mapping $\pi: J_j\to M_j$ and to the partition 
of $M_j$ into $i_j$ intervals $L_i$, we can refine the partition by at most $i_j+(k-1)(i_j-1)$ intervals in 
$M_j$ (but they are also intervals in $[n]$) whose images by $\pi^{-1}$ are intervals in $J_j$ 
(and so in $[n]$). Taking 
all these refinements for $j=1,2,\dots,r$, we get a partition of $[n]$ in at most 
$\sum_{j=1}^r(i_j+(k-1)(i_j-1))<\sum_{j=1}^r ki_j=kl\le kK^*$ intervals whose images by $\pi^{-1}$ are intervals 
in $[n]$ refining the partition $J_1<J_2<\dots<J_r$.
\kduk

\begin{prop}\label{wpo}
For every fixed $k,K\in\N$, the set 
$$
\{\pi\in S:\ \st(\pi)<K\;\&\;s_k(\pi)<K\}
$$
is a wpo with respect to $\prec$.
\end{prop}
\duk
In view of Lemmas~\ref{easy} and \ref{substituce}, it suffices to show that $(H^+_k\cup H^-_k,\prec)$ 
is a wpo. It is enough to show that $(H^+_k,\prec)$ is a wpo. Using $k$-decompositions, we represent each 
$\pi\in H^+_k$ by a word over $\Sigma=\ind^+_1\cup\dots\cup\ind^+_{k-1}$. Now, denoting $\le_s$ the ordering by 
subsequence, it follows from Theorem~\ref{higman} that $(\Sigma^*,\le_s)$ is a wpo and this implies 
that $(H^+_k,\prec)$ is a wpo.
\kduk

\begin{prop}
For every $0<\ep\le 1$, the set $K_{2-\ep}$ is countable.
\end{prop}
\duk
Let an $\ep$, $0<\ep\le 1$, and a CPC $X\in K_{2-\ep}$ be given. It suffices to show that 
the antichain of permutations 
$M(X)$ is finite. We have $|X\cap S_n|<2^{n-1}$ for some $n>1$ and, by Proposition~\ref{struktura}, 
$\st(X)<K$ and 
$s_k(X)<K$ for some constants $k,K\in\N$. By Lemma~\ref{spojitost}, $\st(M(X))<K+2$ and $s_k(M(X))<K+2$. 
By Proposition~\ref{wpo}, $M(X)$ is finite. 
\kduk

\noindent
This finishes the proof of the inequality $\kappa\ge 2$. In fact, we have proved that the set 
$$
\{X: \mbox{$X$ is a CPC such that $|X\cap S_n|<2^{n-1}$ for some $n\in\N$}\}
$$
is countable. It is likely that $K_2$ is countable. 

\section{The upper bound of Theorem~\ref{hlavni}}

Atkinson, Murphy and Ru\v skuc \cite{atki02} introduced an infinite antichain of permutations 
$$
U=\{\mu_7,\mu_9,\mu_{11},\dots\}
$$ 
where
\begin{eqnarray*}
\mu_7&=&4,7,6|1,5,3,2\\
\mu_9&=&6,9,8|4,7|1,5,3,2\\
\mu_{11}&=&8,11,10|6,9,4,7|1,5,3,2\\
 &\vdots&  \\
\mu_{2k+5}&=& 2k+2,2k+5,2k+4|2k,2k+3,2k-2,2k+1,\dots,6,9,4,7|1,5,3,2\\
 &\vdots&
\end{eqnarray*}
The initial segment in $\mu_{2k+5}$ is $2k+2,2k+5,2k+4$, the final segment is $1,5,3,2$, and in the middle segment
the sequences $2k,2k-2,\dots,4$ and $2k+3,2k+1,\dots,7$ are interleaved. (In fact, we have reversed the 
permutations of 
\cite{atki02}). We reprove, using a different argument than in \cite{atki02}, that $\mu_i$ form an antichain. 
We associate with $\pi\in S_n$ a graph $G(\pi)$ on the vertex set 
$\{(i,\pi(i)):\ i\in[n]\}$, in which $(i,\pi(i))$ and $(j,\pi(j))$ are adjacent if and only if 
$i<j$ and $\pi(i)<\pi(j)$. It is clear that $\pi\prec\sigma$ implies $G(\pi)\le_g G(\sigma)$ where $\le_g$ is the 
subgraph relation (this holds even with the induced subgraph relation). A {\em double fork} $F_i$ is the 
tree on $i$ vertices, $i\ge 6$, that is obtained by appending pendant vertex both to the second and to the 
penultimate vertex of a path with $i-2$ vertices. It is easy to see that $(\{F_i:\ i\ge 6\},\le_g)$ is an 
antichain.

\begin{lema}\label{antiret}
$(U,\prec)$ is an antichain. Moreover, 
$$
(\{123, 3214, 2143, 15432\}\cup U,\prec)
$$
is an antichain. 
\end{lema}
\duk
For every $i=7,9,11,\dots$, $G(\mu_i)=F_i$. Since double forks form an antichain to $\le_g$, so do 
the permutations $\mu_i$ to $\prec$. It is clear that
the four new short permutations form an antichain and none contains any $\mu_i$. 
$G(123)$ is a triangle, $G(2143)$ is a quadrangle and $G(15432)$ has a vertex of degree $4$, 
and therefore none of the three permutations is contained in any $\mu_i$. That $3214\not\prec\mu_i$ for 
every $i$ is easily checked directly. 
\kduk

\begin{prop}\label{enumerace}
Let $s_n=|S_n(123, 3214, 2143, 15432)|$. Then
$$
\sum_{n\ge 1}s_nx^n=\frac{x^5+x^4+x^3+x^2+x}{1-x-2x^2-2x^3-x^4-x^5}.
$$
As $n\to\infty$, $s_n\sim c(2.33529\dots)^n$ where $c>0$ is a constant and $2.33529\dots$ is the only real root 
of $x^5-x^4-2x^3-2x^2-x-1$.
\end{prop}
\duk
We denote $S^*_n=S_n(123, 3214, 2143, 15432)$ and partition $S_n^*$ in five sets $A_n,\dots, E_n$ as 
follows. For $n\ge 2$ and $\pi\in S^*_n$, we let $\pi\in A_n\Leftrightarrow\pi(1)=n-1$, $\pi\in B_n\Leftrightarrow\pi(1)=n-2$, 
$\pi\in C_n\Leftrightarrow\pi(1)\le n-3$, $\pi\in D_n\Leftrightarrow\pi(1)=n\;\&\;\pi(2)\ge n-3$, and 
$\pi\in E_n\Leftrightarrow \pi(1)=n\;\&\;\pi(2)\le n-4$. We denote $|A_n|=a_n,\dots,|E_n|=e_n$.  
Notice that for every $n\in\N$ and $\pi\in S^*_n$, $\pi^{-1}(n)\le 3$. For if $\pi^{-1}(n)\ge 4$, 
the first three values of $\pi$ have an ascend or all are descending, and $123\prec\pi$ or $3214\prec\pi$. 
Thus every $\sigma\in S^*_{n+1}$ 
arises from some $\pi\in S^*_n$ by inserting the value $n+1$ on one of the three sites: in front of the whole $\pi$ 
(site 1), between the first two values of $\pi$ (site 2) or between the second and the third value of $\pi$ 
(site 3). We discuss the cases depending on in which set $\pi$ lies.

In all five cases we can insert $n+1$ on site 1. With the exception of the case $\pi\in D_n$, we cannot insert 
$n+1$ on site 3 because this would give $123\prec\sigma$ or $2143\prec\sigma$ or $15432\prec\sigma$. If 
$\pi\in C_n$, we cannot insert $n+1$ on site 2 because this would give $123\prec\sigma$ or $15432\prec\sigma$.
One can check that there are no other restrictions on the insertion of $n+1$. 
Hence $\pi\in A_n$ produces 
two $\sigma$'s, one in $D_{n+1}$ and the other in $B_{n+1}$; $\pi\in B_n$ produces also
two $\sigma$'s, one in $D_{n+1}$ and the other in $C_{n+1}$; $\pi\in C_n$ produces one $\sigma$ in $E_{n+1}$; 
$\pi\in D_n$ produces three $\sigma$'s, one in $D_{n+1}$ and two in $A_{n+1}$; and $\pi\in E_n$ produces 
two $\sigma$'s, one in $D_{n+1}$ and the other in $A_{n+1}$. From this we obtain the recurrences 
$a_{n+1}=2d_n+e_n$, $b_{n+1}=a_n$, $c_{n+1}=b_n$, $d_{n+1}=a_n+b_n+d_n+e_n$, and $e_{n+1}=c_n$. 

We set 
$(a_1,b_1,c_1,d_1,e_1)=(0,0,0,0,1)$, which gives correctly $(a_2,b_2,c_2,d_2,e_2)=(1,0,0,1,0)$. Let 
$v=(0,0,0,0,x)$ be the vector of initial conditions for $n=1$ and $M$ be the $5\times 5$ transfer matrix
$$
M=\left(
\begin{array}{rrrrr}
0&0&0&2x&x\\
x&0&0&0&0\\
0&x&0&0&0\\
x&x&0&x&x\\
0&0&x&0&0
\end{array}
\right).
$$
For the generating functions $A=\sum_{n\ge 1}a_nx^n,\dots, E=\sum_{n\ge 1}e_nx^n$, the recurrences give 
relation
$$
(A,B,C,D,E)^T=(I+M+M^2+\cdots)v^T=(I-M)^{-1}v^T. 
$$
From this, since $s_n=a_n+b_n+c_n+d_n+e_n$, 
\begin{eqnarray*}
\sum_{n\ge 1}s_nx^n&=&A+B+C+D+E=(1,1,1,1,1)(I-M)^{-1}v^T\\
&=&\frac{x(x^4+x^3+x^2+x+1)}{1-x-2x^2-2x^3-x^4-x^5}.
\end{eqnarray*}

One can check that $2.33529\dots$ is the dominant root of the reciprocal polynomial $x^5-x^4-2x^3-2x^2-x-1$ 
of the denominator. The asymptotics of $s_n$ follows from the standard facts on asymptotics of coefficients 
of rational functions. 
\kduk

\noindent
We obtain the recurrence $s_1=1$, $s_2=2$, $s_3=5$, $s_4=12$, $s_5=28$, and 
$s_n=s_{n-1}+2s_{n-2}+2s_{n-3}+s_{n-4}+s_{n-5}$ for $n\ge 6$.
The first values of $s_n$ are:
$$
(s_n)_{n\ge 1}=(1,2,5,12,28,65,152,355,829,1936,4521,10558,\dots).
$$

\begin{prop}
For every $\ep>0$, the set $K_{2.33529\dots+\ep}$ is uncountable.
\end{prop}
\duk
The set of CPC's
$$
\{S(\{123, 3214, 2143, 15432\}\cup V):\ V\subset U\}
$$
is uncountable, due to Lemma~\ref{antiret} and the 1-1 correspondence between CPC's and antichains of permutations, 
and 
$$
|S_n(\{123, 3214, 2143, 15432\}\cup V)|\le |S_n(123, 3214, 2143, 15432)|=s_n.
$$
By Proposition~\ref{enumerace} we know that for any $\ep>0$, $s_n<(2.33529\dots+\ep)^n$ for every $n>n_0$. 
\kduk

\noindent
Thus $\kappa\le 2.33529\dots$ and the proof of Theorem~\ref{hlavni} is complete. More restrictions can be added 
to the $\{123, 3214, 2143, 15432\}$-avoidance and the bound $\kappa\le 2.33529\dots$ can be 
almost surely improved but 
the question is by how much. It seems not very likely that one could prove this way that $\kappa\le 2$. 

We conclude with some comments on our choice of the four permutations $123, 3214, 2143$, and $15432$. 
By the results in \cite{atki02}, if $(S(\pi,\rho),\prec)$ is not a wpo, where $\pi\in S_3$, $\rho\in S_4$ and 
$\pi\not\prec\rho$, then $(\pi,\rho)$ equals, up to obvious symmetries, to $(123,3214)$ or 
$(123,2143)$. In \cite{atki02} it is also observed that $S(123, 3214, 2143)\supset U$ and so 
$(S(123, 3214, 2143),\prec)$ is not a wpo. We have employed one more restriction: From the 28 permutations in 
$S_5(123, 3214, 2143)$, only $15432$ is not contained in infinitely many $\mu_i$. The enumeration 
$|S_n(123)|=C_n$, where $C_n$ is the $n$-th Catalan number, is a classic result (see Stanley~\cite{stan}); $C_n$
have exponential growth $4^n$. West \cite{west96} proved that $|S_n(123,3214)|=|S_n(123,2143)|=F_{2n}$ where 
$(F_n)_{n\ge 1}=(0,1,1,2,3,5,8,13,\dots)$ are Fibonacci numbers. $F_{2n}$ grow as 
$((3+\sqrt{5})/2)^n=(2.61803\dots)^n$. Using simpler arguments than those in the proof of 
Proposition~\ref{enumerace}, we can prove that the numbers $t_n=|S_n(123, 3214, 2143)|$ follow the recurrence
$t_1=1$, $t_2=2$ and $t_n=2t_{n-1}+t_{n-2}$ for $n\ge 3$. Thus $t_n$ grow as $(1+\sqrt{2})^n=(2.41421\dots)^n$.

\end{document}